\documentclass[11pt,letterpaper]{article}
\usepackage{ifthen,latexsym,amssymb,amsmath}
\usepackage{graphicx}
\usepackage{booktabs}
\usepackage{url}

\newcommand{\QQ}{\mathcal{Q}}
\newcommand{\LL}{\mathcal{L}}

\newif\ifnotesw\noteswtrue

\noteswfalse   

\newcommand{\beq}[1]{\begin{equation}\label{eq:#1}}
\newcommand{\eeq}{\end{equation}}

\newtheorem{theorem}{Theorem}
\newcommand{\bth}[2][nothing]{\ifthenelse{\equal{#1}{nothing}}
 {\begin{theorem}} {\begin{theorem}[#1]}\label{th:#2}}

\newtheorem{lemma}[theorem]{Lemma}
\newcommand{\blm}[2][nothing]{\ifthenelse{\equal{#1}{nothing}}
 {\begin{lemma}} {\begin{lemma}[#1]}\label{lm:#2}}

\newtheorem{problem}[theorem]{Problem}
\newcommand{\bpr}[2][nothing]{\ifthenelse{\equal{#1}{nothing}}
 {\begin{problem}} {\begin{problem}[#1]}\label{pr:#2}}

\newcommand{\qed}{\nolinebreak\mbox{\hspace{5 true pt}%
  \rule[-0.85 true pt]{3.9 true pt}{8.1 true pt}}}

\newcommand{\nothree}{no-3-in-a-line }

\title{Martin Gardner's minimum no-3-in-a-line problem}
\author{Alec S. Cooper, Oleg Pikhurko, John R. Schmitt, Gregory S. Warrington}
\begin{document}
\date{}

\maketitle
\begin{abstract}
  In Martin Gardner's October 1976 Mathematical Games column in {\it
    Scientific American}, he posed the following problem: ``What is
  the smallest number of [queens] you can put on an [$n \times n$ chessboard]
  such that no [queen] can be added without creating three in a row, a
  column, or a diagonal?''  We use the Combinatorial Nullstellensatz
  to prove that this number is at least $n$, except in the case when $n$ is congruent to $3$ modulo $4$, in which case one less may suffice.  A second, more
  elementary proof is also offered in the case that $n$ is even.
\end{abstract}

\section{Introduction.}

In Martin Gardner's October 1976 Mathematical Games column in {\it
  Scientific American}, he introduced this combinatorial chessboard
problem: What is the minimum number of counters that can be placed on
an $n\times n$ chessboard, no three in a line, such that adding one
more counter on any vacant square will produce three in a line?  He
dubbed the problem the {\it minimum \nothree problem}.  

\begin{figure}[htbp]
  \centering 
  {\includegraphics[height=150px]{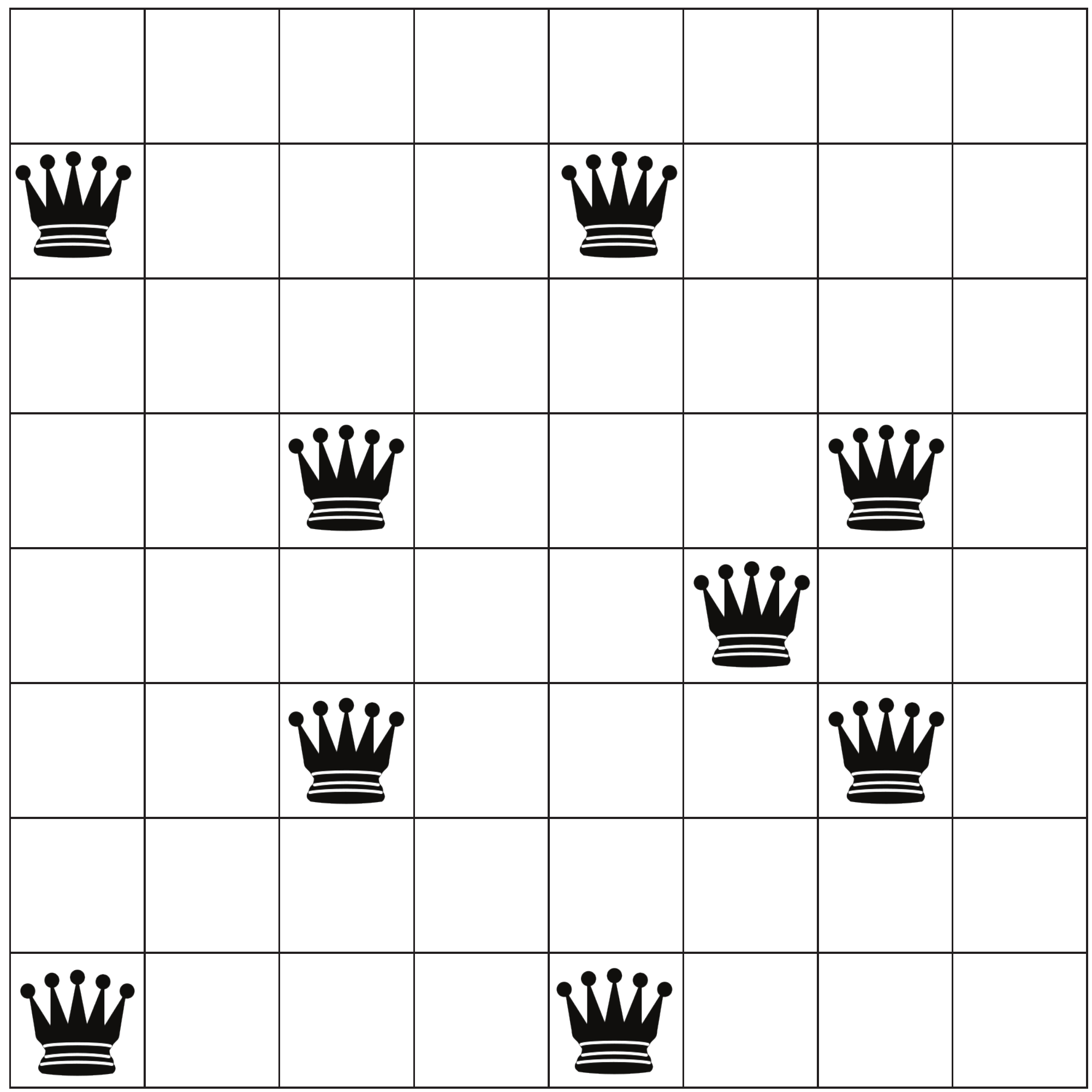}}\qquad
  {\includegraphics[height=150px]{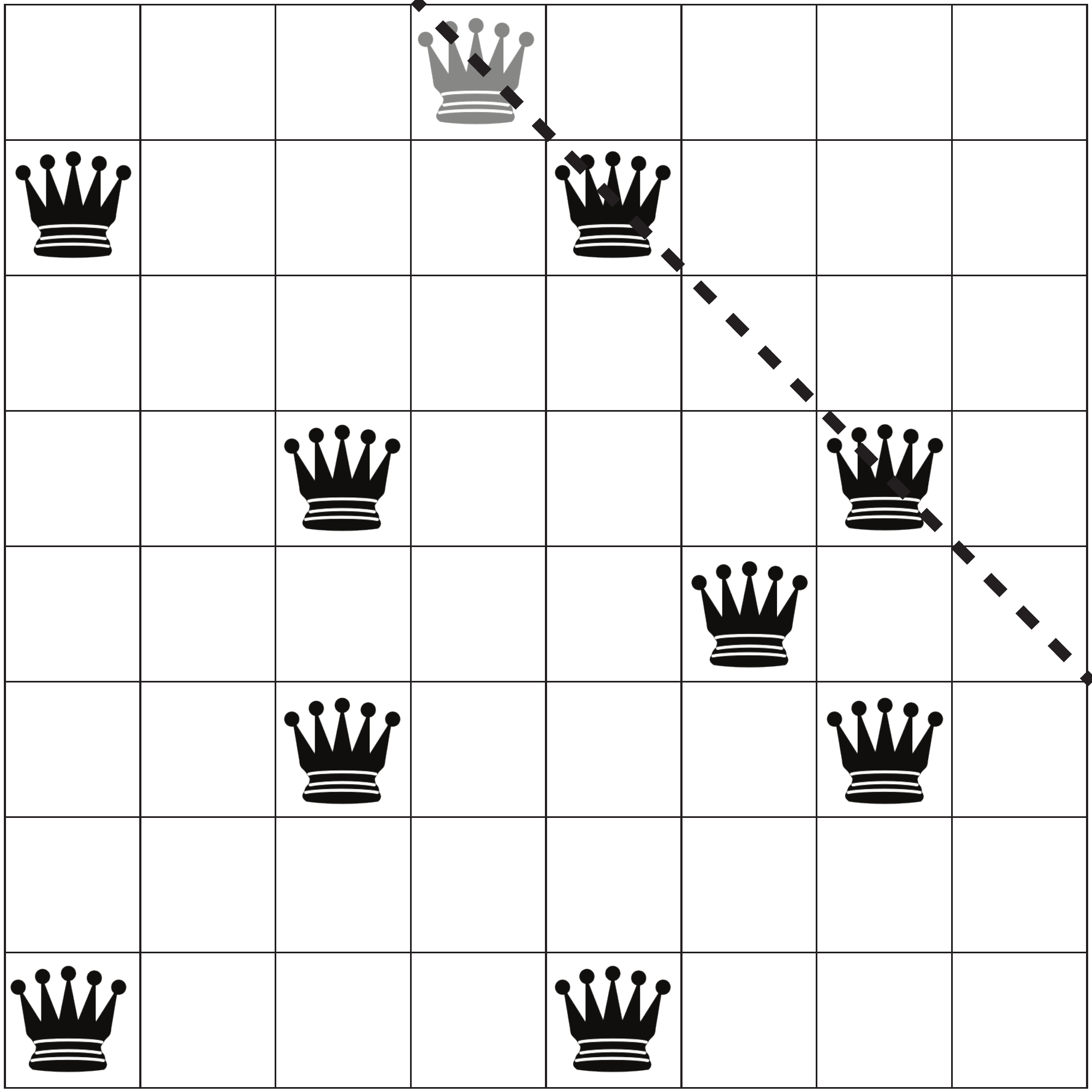}}
  \caption{Maximal $8 \times 8$ placement}
  \label{fig:8by8ex}
\end{figure}

Figure~\ref{fig:8by8ex} shows an $8\times 8$ chessboard with an
initial placement of $9$ black queens with no three in a line.  This
placement is maximal, that is, any additional queen will create three
in a line.  The figure illustrates the corresponding `three-in-a-line'
 created when
an additional queen, shown in a distinct shading, is placed in the fourth
column and eighth row.  This particular placement is also of minimum size (where {\it size} of a placement is the number of queens in the placement), that is, there 
is no placement with eight or fewer queens meeting the requirements.

Gardner makes the following observation~\cite[Chapter 5,
pg. 71]{gardner:klondike}:
\begin{quotation}
  If `line' is taken in the broadest sense --- a straight line of any
  orientation --- the problem is difficult\ldots The problem is also
  unsolved if `line' is restricted to orthogonals and diagonals.
\end{quotation}

In this paper we provide a lower bound for this latter \emph{queens
  version} of the problem.

\begin{theorem}\label{theorem:main}
  For $n \geq 1$, the answer to Gardner's \nothree problem is at least $n$, except in the case when $n$ is congruent to $3$ modulo $4$, in which case one less may suffice.
\end{theorem}

We offer an elementary, ad hoc proof in the case of $n$ even (the
approach yields a lower bound of only $n-1$ when $n$ is odd).  This
proof, similar to an incomplete argument of John Harris~\cite{harris-letter},
is provided in Section~\ref{section:Harrisproof}.

The proof of Theorem~\ref{theorem:main} for arbitrary $n$ ultimately
relies on the Nullstellensatz.  Hilbert's ``zero-locus
theorem''~\cite{hilbert:93} is a foundational result that connects
geometry and algebra.  In~\cite{alon:99}, Alon leverages a special
case of Hilbert's theorem to prove a \emph{Combinatorial
  Nullstellensatz} (reproduced here as Theorem~\ref{theorem:CN}) that
is ideally suited for obtaining lower bounds on restricted-sum sets
and other similar objects (see \cite[Chapter 9]{tao+vu:10}).

The proof we present in Section~\ref{section:CNproof} using the
Combinatorial Nullstellensatz is inspired by a similar proof of Alon and F\"uredi~\cite{alon-furedi:93}; their proof gives a result
about the number of hyperplanes needed to cover all but one of the
vertices of the hypercube (see \cite[Theorem~6.3]{alon:99}).  We
believe that our proof serves as a nice illustrative application of the
Combinatorial Nullstellensatz.

We may arrive at lower bounds that are weaker than those promised by
Theorem \ref{theorem:main} quite quickly.  If we make the observation that each of the $q$ queens
`covers' at most $4n-4$ squares and each of the $n^2$ squares requires
either two queens to `cover' it or one queen to occupy it, a lower
bound of $\frac{n}{2}$ follows~\cite{dan}.  This last observation can
be strengthened by noting that only a few queens can cover $4n-4$
squares.  However, any queen covers at `worst' $3n-3$ squares, though still we could not push this line of argument to get us to $n$.

Prior to our proof of Theorem~\ref{theorem:main}, we discuss some history of the problem
drawn from Gardner's notes and correspondence pertaining to his
writing of the {\em Scientific American} column
\cite{gardner:archives} -- this is done in Section~\ref{section:history}.  Section~\ref{section:CNproof}, as mentioned
above, presents a proof of Theorem~\ref{theorem:main} using the
Combinatorial Nullstellensatz.  We also offer a more elementary proof
in Section~\ref{section:Harrisproof}.

\section{History.}
\label{section:history}

Gardner and some of his readers found good placements -- a placement is {\it good} if it does not contain three queens in a line and loses this property upon the addition of a queen to an unoccupied square -- via
pencil-and-paper; others conducted computer searches.  We also
conducted computer searches, though with computing power that is
better than it was 35 years ago.  Collectively, these results are
contained in Table \ref{table:1}; $m_3(n)$ denotes the answer to
Gardner's \nothree problem on an $n\times n$ chessboard.  A
bold-faced entry in the second row indicates that an improvement was
made to previous knowledge.

Theorem~\ref{theorem:main} is not ``tight'' for small values of $n$.
The data suggest for $n$ odd and $n \geq 3$ that we should have
$m_3(n) \geq n+1$.  Our search for good placements was done via
brute-force search\footnote{The C code that performed this search is
  available in the source package for~\cite{arx} at
  \url{http://arXiv.org}.}.  As such, and to illustrate the
computational challenges involved, our program took around 900 3GHz-CPU hours
to confirm that there is no good placement of $11$ queens on an $11
\times 11$ chessboard.  We estimate that the
corresponding search for a $13 \times 13$ chessboard using our program
would require at least 70 thousand 3GHz-CPU hours.

\begin{table}
  \begin{center}
    \begin{tabular}{cccccccccc}\toprule
      $n$       & 1 & 2 & 3 & 4 & 5 & 6 & 7 & 8 & 9 \\
      $m_3(n)$  & 1 & 4 & 4 & 4 & 6 & 6 & 8 & {\bf 9} & {\bf 10}\\\midrule
      $n$ & 10 &       11 & 12 & 13      & 14 & 15 & 16 & 17 & 18\\
      $m_3(n)$ & {\bf 10} & {\bf 12} & 12 & {\bf [13,14]} & {\bf [14,16]} & 
      {\bf [14,16]} & {\bf [16,18]} & {\bf [17,20]} & {\bf [18,20]}\\\bottomrule
    \end{tabular}
  \end{center}
  \caption{$m_3(n)$, for small values of $n$.  Brackets indicate lower and upper bounds.}
  \label{table:1}
\end{table}

\begin{figure}[htbp]
  \centering 
{\includegraphics[height=130px]{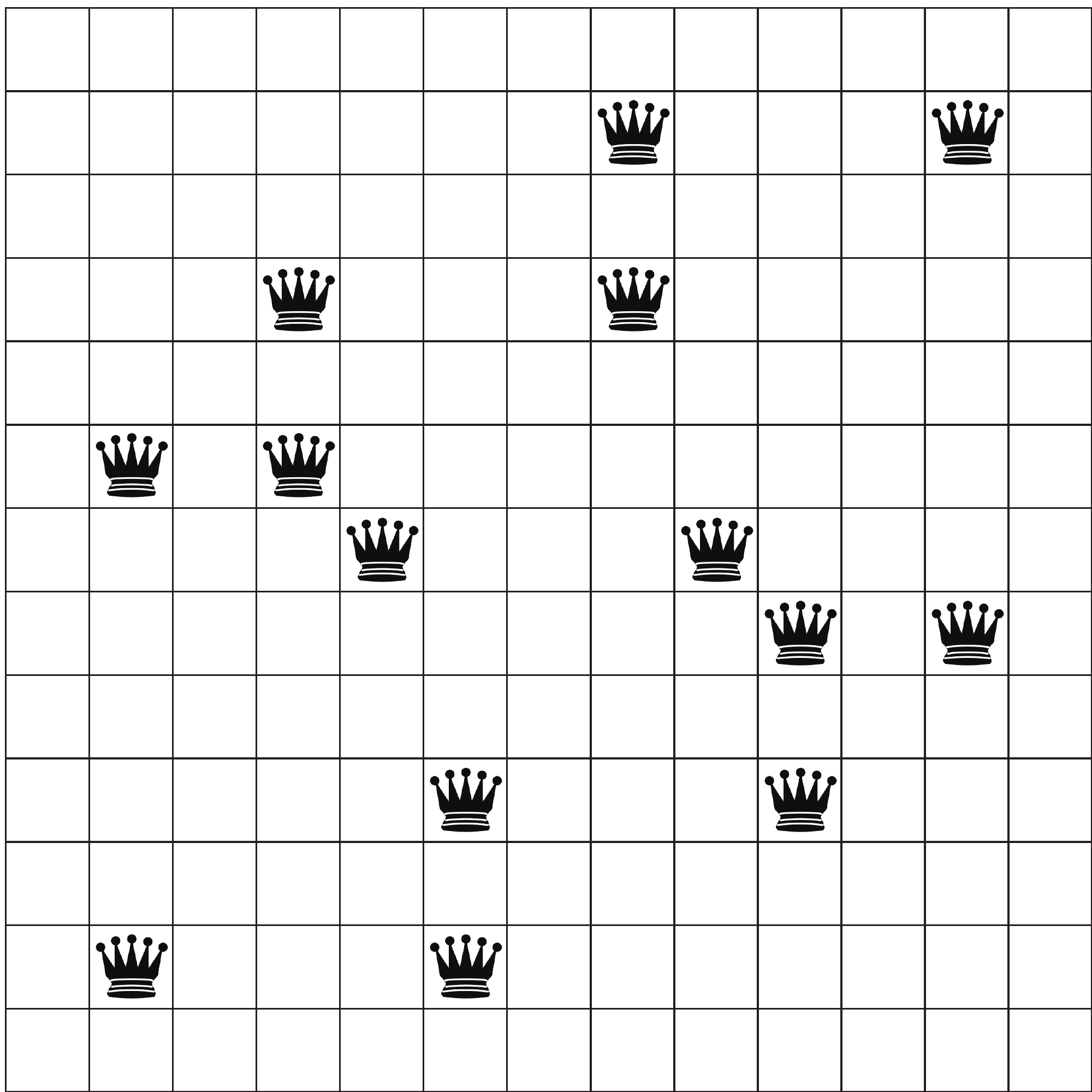}}\quad
{\includegraphics[height=140px]{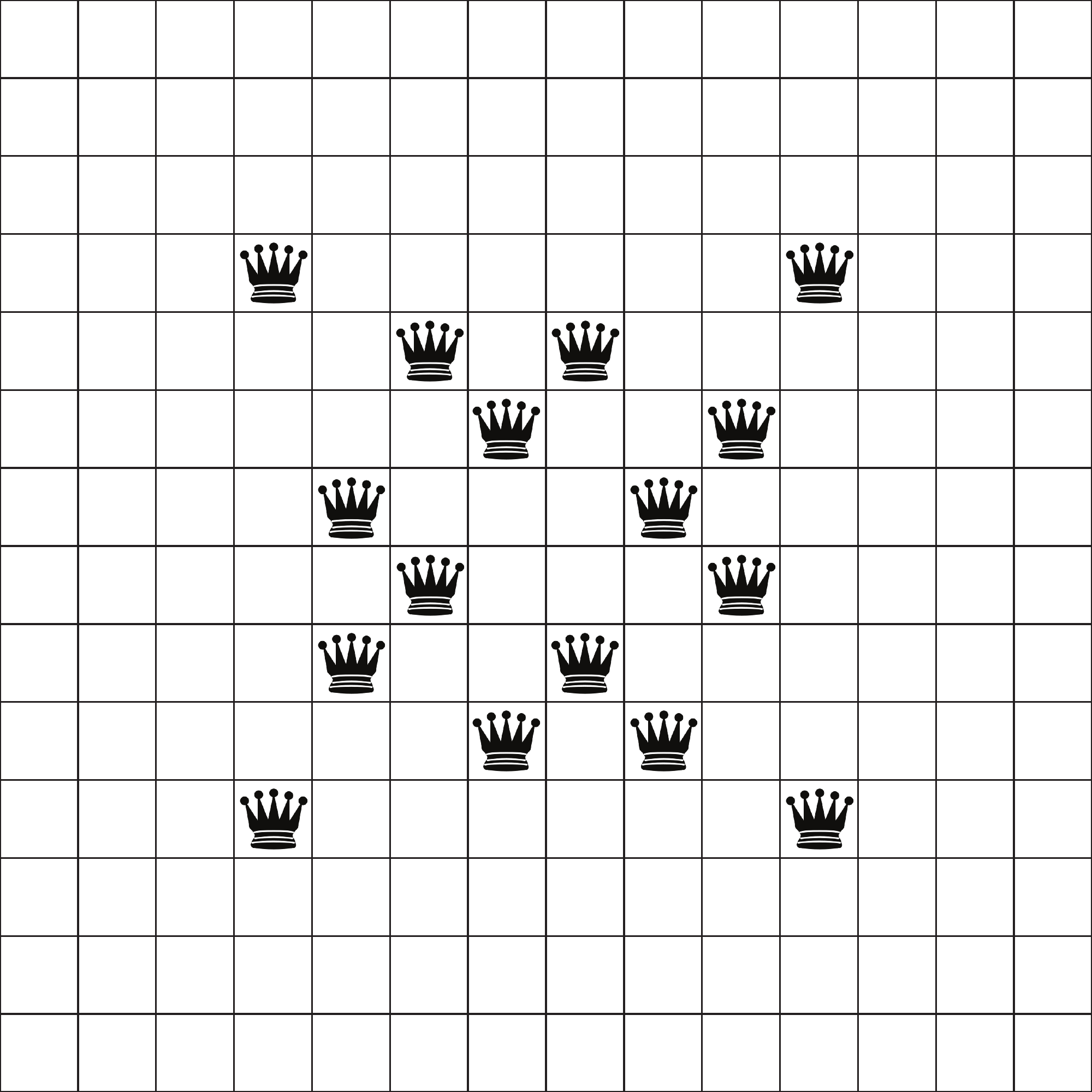}}\quad
    {\includegraphics[height=150px]{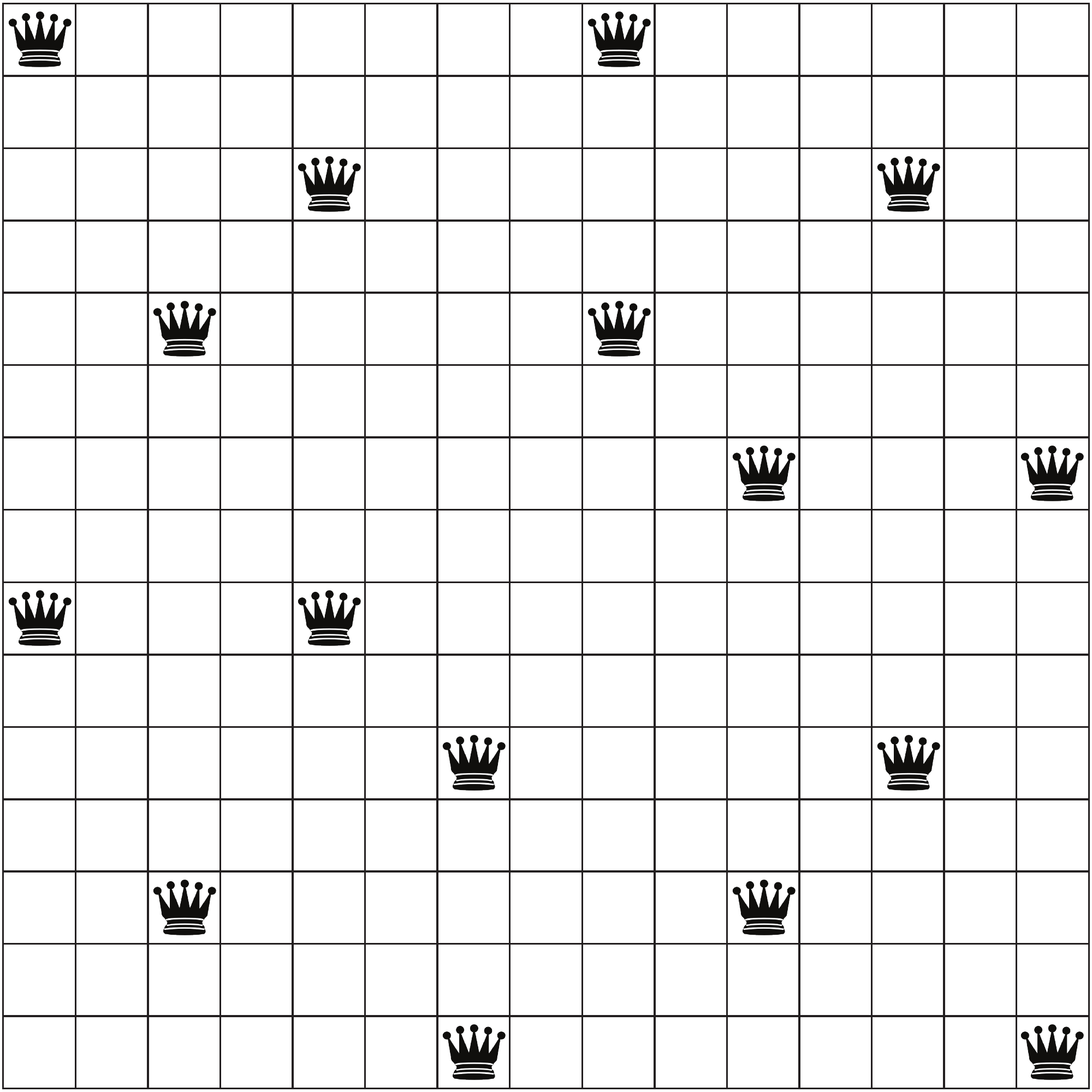}}
     \caption{Maximal placements: $14$ queens for $n=13$; $16$ queens for
       $n\in\{14,15\}$.}
      \label{fig:goodex}
\end{figure}

\begin{figure}[htbp]
  \centering 
 {\includegraphics[height=130px]{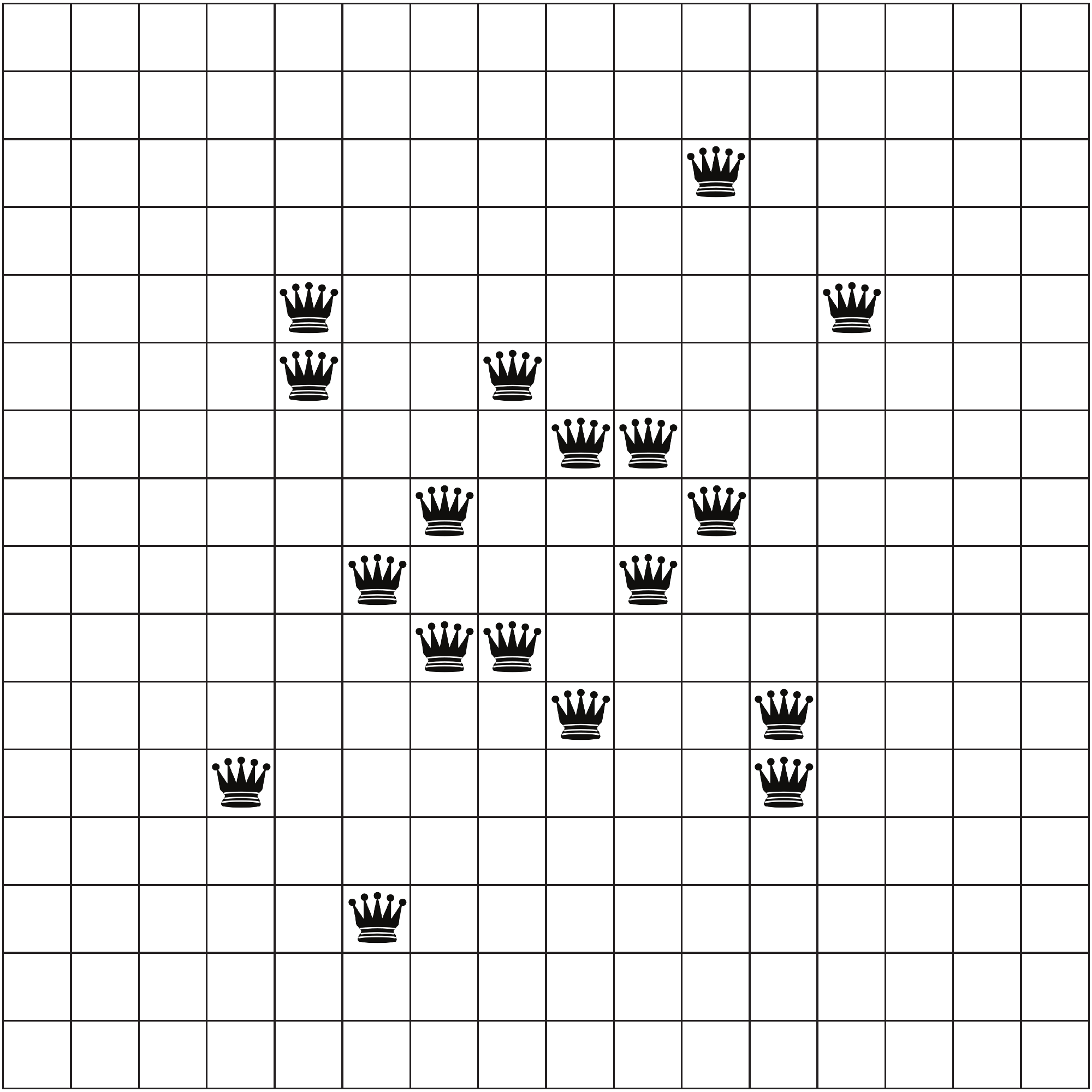}}\quad
{\includegraphics[height=140px]{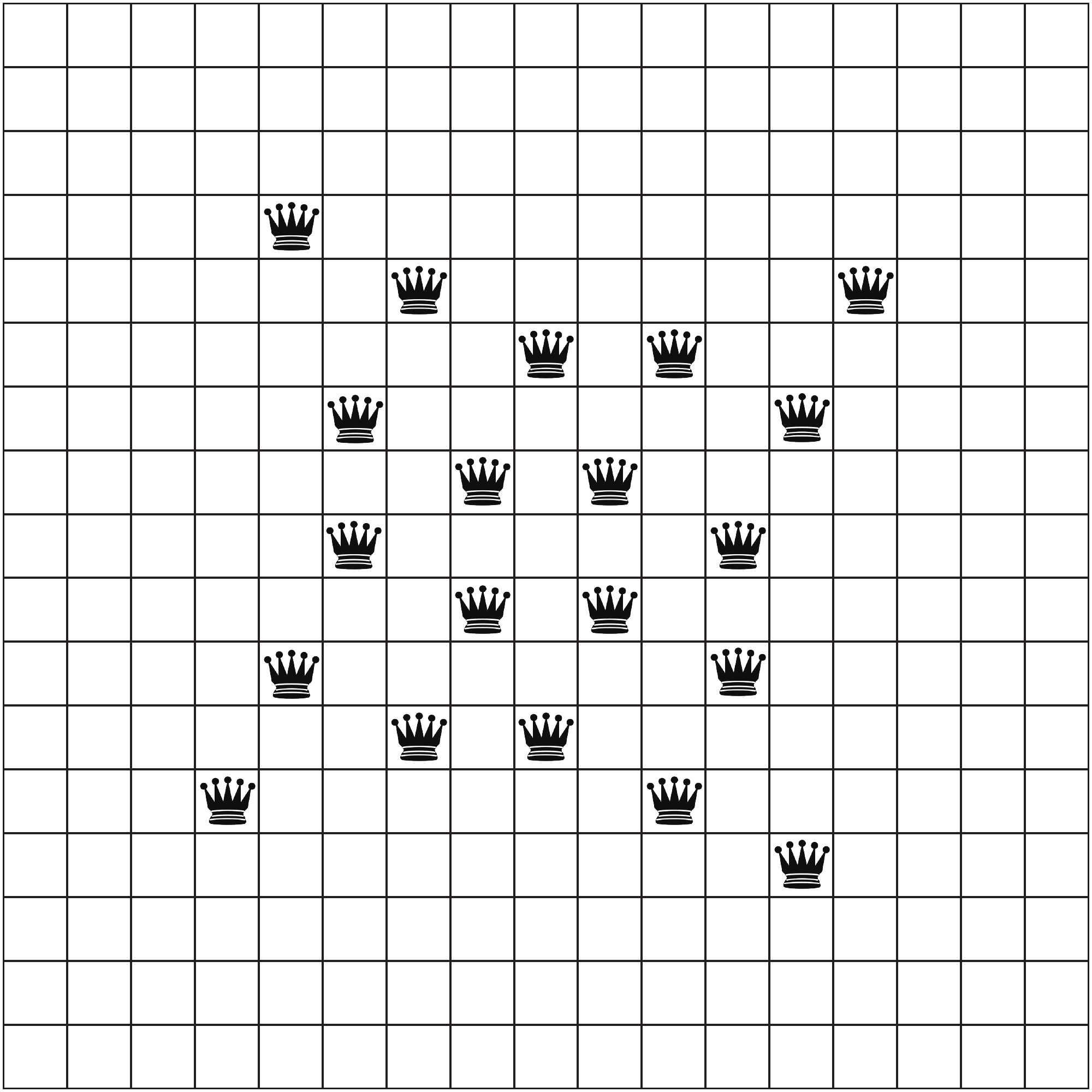}}\quad
   {\includegraphics[height=150px]{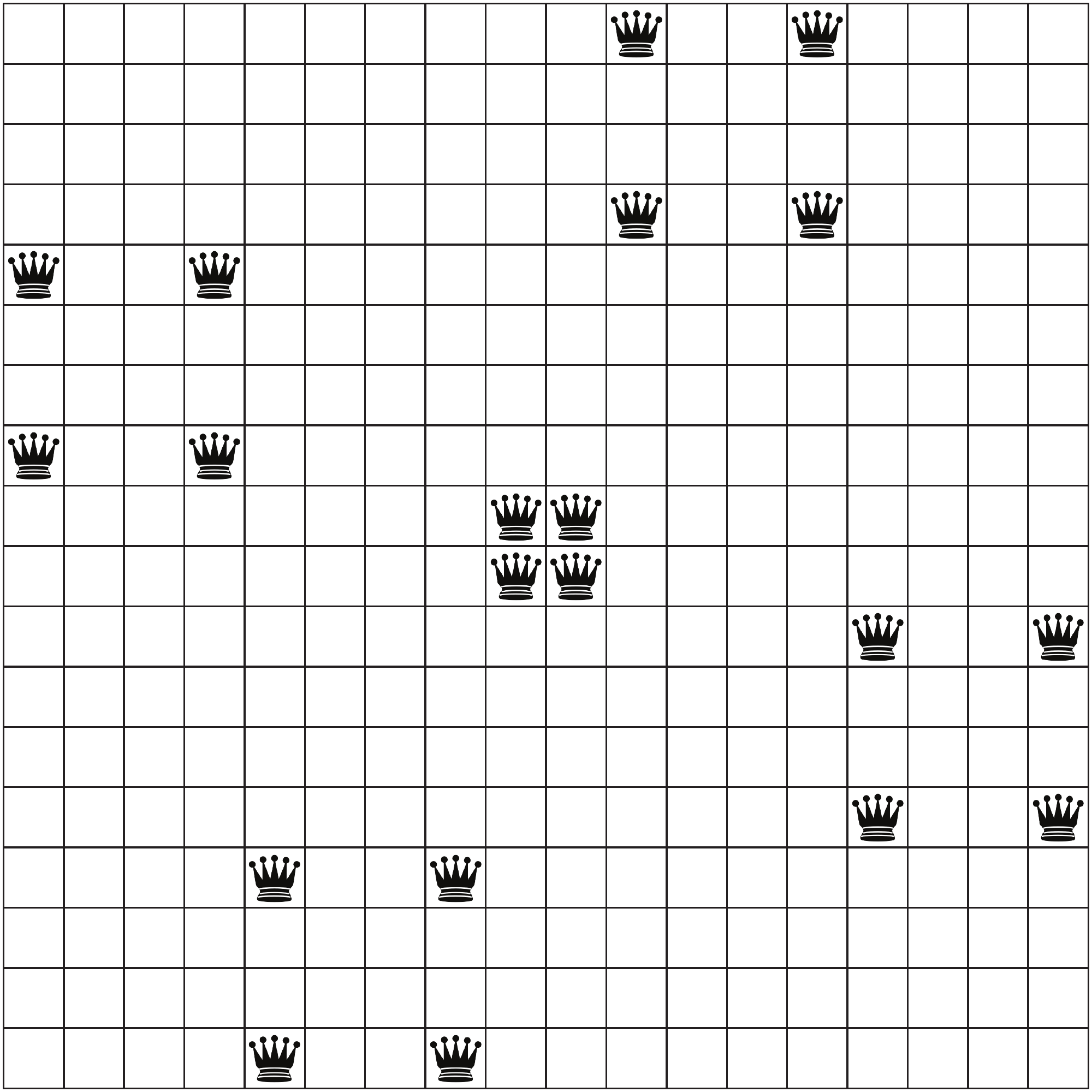}}
     \caption{Maximal placements: $18$ queens for $n=16$; $20$ queens for
       $n\in\{17,18\}$.}
      \label{fig:goodex-ii}
\end{figure}

In that October column (and in an addendum \cite{gardner:klondike} to it), Gardner gave a few
results on the queens version.  These included placements of queens on
chessboards  $3 \times 3$ through $12 \times 12$, which provided upper bounds on this
number.  His archives also contain a good placement of $52$ queens on
a $48 \times 48$ chessboard~\cite{gardner:archives}.  Gardner also stated that
John Harris of Santa Barbara, CA (who, we later learned, was a
frequent correspondent of Gardner's) was able to show that the minimum number of queens needed for an $n \times n$ chessboard is at least $n$, except
when $n$ is congruent to 3 modulo 4, in which case
it could be one less.  Gardner did not supply Harris' argument.  These
results were the ``jumping off'' point of our investigations ---
mostly, we wondered what Harris' argument was.

Subsequent to obtaining our results that confirm and improve upon
those of Harris, we were able to obtain copies of Gardner's notes and
correspondences concerning this problem \cite{gardner:archives}.
These are a small fraction of the 60 linear feet(!) of notes and
correspondences archived at Stanford University that pertain to his
writing of the Mathematical Games column.  (These materials were a
gift to Stanford by Gardner in 2002.)  There are numerous
carbon-copies of letters that Gardner wrote to other mathematicians,
as well as readers, and copies of letters they wrote to him about the
problem.  Chronologically first is a letter, dated June 2, 1975, that
Gardner wrote to the world-renowned John H. Conway.  In it, he states
that the problem occurred to him while considering a game of the
mathematician Stanislaw Ulam (though not the game that commonly goes
by the name Ulam's game) --- the game had appeared in an earlier
column.  The game consists of taking turns ``putting a counter on an
$n \times n$ [chessboard] until one person wins by getting 3 in line,
orthogonally or diagonally.''  In the weeks that followed were letters
to and from Bill Sands (then a Ph.D. student at U. Manitoba, now at
U. Calgary), who independently suggested the problem, and John Harris,
including one that sketches some ideas for the above-mentioned claim.
Subsequent letters from readers (that appeared after the October 1976
column) contained their best solutions to the problem for small chessboards; some of Gardner's notes do the same.

Of course, the reader may be more aware of some related or similarly
worded problems.  Gardner mentioned one of them in that month's
column, the {\it maximum \nothree problem}, that is, what is the
maximum number of counters (or queens) one can place on an $n \times n$ chessboard so that there are no three in a line?  Here an easy upper
bound of $2n$ follows from the pigeonhole principle as each of the $n$
columns may contain at most $2$ counters --- Guy and Kelly
\cite{guy+kelly:68} showed that one is `unlikely' to find any with more
than $\sim1.87n$ queens -- this was later corrected to $\sim1.81n$ queens  (see  \cite[{\fontfamily{ppl}\selectfont A000769}]{sloane}).  Another related problem is the {\it queens
  domination problem}.  In this problem, one asks for the minimum
number of queens needed so that each square of the chessboard is
either occupied or attacked.  There are two versions of this problem,
one where the queens are non-attacking and the other where this
restriction is lifted, see \cite{burger+mynhardt:03} and
\cite{fernau:10} for some results.

%

\section{Proof of Main Theorem via the Combinatorial
  Nullstellensatz.}\label{section:CNproof}

In this section we prove Theorem~\ref{theorem:main} using the
Combinatorial Nullstellensatz.  To begin, we give a brief discussion of
the theorem to be applied and its statement.

The Fundamental Theorem of Algebra tells us that a degree-$t$
polynomial $f(x)$ contained in a polynomial ring $F[x]$ has at most
$t$ zeros.  Said another way, for any set $S$ contained in $F$ of
cardinality greater than $t$, there is an element $s \in S$ such that
$f(s)$ is nonzero.  One may think of this as saying, either a
polynomial is zero everywhere or it is zero in very few places.  The
following theorem, known as the Combinatorial Nullstellensatz,
generalizes this fact to polynomials of several variables --- it is
due to Alon \cite[Theorem~1.2]{alon:99}.  We may think of it as saying
that a multivariable polynomial that isn't zero everywhere has a non-root in a box of large
enough volume.

\begin{theorem}\label{theorem:CN}[Combinatorial Nullstellensatz, Theorem 1.2~\cite{alon:99}]
  Let $F$ be an arbitrary field, and let $f = f ( x_1 , \ldots , x_n
  )$ be a polynomial in $F [ x_1 , \ldots , x_n ]$.  Suppose the
  degree $\deg (f)$ of $f$ is $\sum_{i = 1}^n t_i$, where each $t_i$
  is a nonnegative integer, and suppose the coefficient of $\prod_{i =
    1}^n x_i^{t_i}$ in $f$ is nonzero.  Then, if $S_1 , \ldots , S_n$
  are subsets of $F$ with $|S_i| > t_i$, there are $s_1 \in S_1 ,
  \ldots , s_n \in S_n$ so that $f ( s_1 , \ldots , s_n ) \neq 0$.
\end{theorem}

So that we might precisely state our results, we introduce some
definitions and notation.  We consider the infinite square ${\mathbb
  Z}$-lattice as a {\it chessboard} and its vertices as {\it squares} of
the chessboard.  A {\it board} $B$ is a finite subset of the
chessboard.  Let $B_n$ denote the board $[1,n] \times [1,n]$.  As we
are interested in the queens version of the problem, the lines that we
concern ourselves with have slope $0, +1, -1,$ or $\infty$ and contain
vertices of the lattice --- so, throughout we use {\it line} to refer
to a line of this type.  Any subset $S$ of the infinite square lattice
may be considered a {\it placement of queens}, or {\it placement} for
short, by imagining a queen on each corresponding square of the
chessboard.  The {\it size} of a placement $S$ is its cardinality
$|S|$.  We say that two queens of a placement ${\cal Q}$ {\it define}
a line if they lie on the same row, column or diagonal.  In such a
way, the placement ${\cal Q}$ {\it defines} a set of lines, the set of
lines defined by the pairs of ${\cal Q}$.  Lastly, we call a placement
{\it good} if does not contain $3$ queens in a line and loses this
property upon the addition of a queen to an unoccupied square.

Let $m_k(n)$ denote the minimum size of a placement on $B_n$ such that
there are no $k$ queens in a line and the placement loses this
property upon the addition of a queen to an unoccupied square of
$B_n$.  As indicated by our title, our focus is on $k=3$.

\begin{figure}[htbp]
  \centering 
     {\scalebox{.15}{\includegraphics{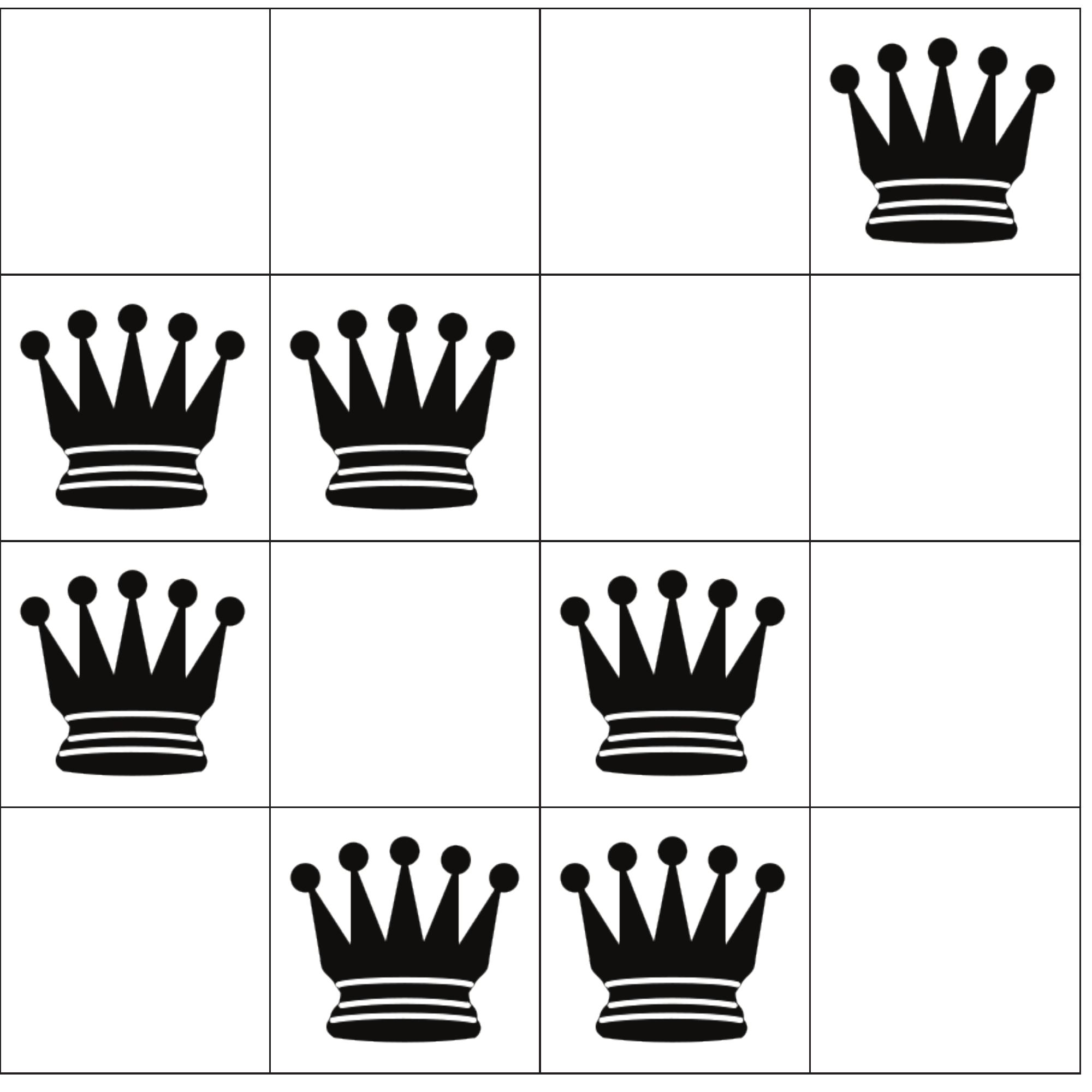}}}
      \caption{Good placement with one queen not collinear with any other.}
      \label{fig:qprimeex}
\end{figure}

We warn the reader that the placements we seek need not have each
queen of a placement on a line with another queen.  See
Figure~\ref{fig:qprimeex} for an example with $n=4$.

We first prove the result for $n = 4k+1$, where $k$ is a nonnegative
integer, as in this case the presentation is cleanest.  We next
establish the result for $n=4k$, and we omit the details for the other
two cases as these are similar.

{\sc Proof of Theorem~\ref{theorem:main}}

Let $n=4k+1$, where $k$ is a positive integer.  (The result is obvious
for $k=0$, i.e., $n=1$.)  Let ${\cal Q}$ be a good placement on $B_n$
with size $q = |{\cal Q}| \leq 4k$.  Our proof will proceed by
constructing a polynomial $f(x,y)$ of total degree $8k$ that vanishes
on each square $(x,y)\in B_n$.  We will then obtain a contradiction
through a suitable application of the Combinatorial Nullstellensatz.

We shall construct $f$ as a product of linear factors of three
different types.  The first type consists of the set of lines defined
by $\QQ$.  Since the placement $\QQ$ is good, every unoccupied square
of $B_n$ is in the zero locus of at least one line of the first type.

As shown in Figure~\ref{fig:qprimeex}, there may be some queens in
$\QQ$ not on any defining line.  Let $\QQ' = \{Q_1,\ldots,Q_{q'}\}$
denote the (possibly empty) subset of queens not collinear with any
other queen in ${\cal Q}$.  For each $Q_i \in \QQ'$ we define a new
line that passes through the square occupied by $Q_i$.  While we are
free to choose any one of the four possible slopes for each line, it is most convenient to distribute the slopes as evenly as
possible.  Hence we choose the slope of the $i^{th}$ line to be 
the $j^{th}$ element of $(0,+1,-1,\infty )$, where $j \equiv i \mod
4$.  Every occupied square is in the zero locus of at least one line of
either of the first two types.

For each of the four possible slopes there are at most $\left\lfloor
  \frac{4k-q'}{2} \right\rfloor$ lines of that slope of the first type
and at most $\left\lceil \frac{q'}{4} \right\rceil$ lines of that
slope of the second type.  These quantities sum to at most $2k$.  As
necessary, define new, distinct lines (of the third type) in each
of the four directions so that there are exactly $2k$ lines of each
slope among the three types.  (The lines of the third type serve only
to facilitate the application of the Combinatorial Nullstellensatz; it
is immaterial which squares they vanish on.)

Let $\LL = \{L_1,\ldots, L_{8k}\}$ be our set of $8k$ lines and let
$l_i=0$ be the equation in variables $x$ and $y$ defining $L_i$.  We
then define
\begin{equation}
  f(x,y)  =  \prod_{i=1}^{8k} l_i \in \mathbb{R}[x,y].
\end{equation}
As desired, the polynomial $f(x,y) = 0$ for every $(x,y)\in B_n$ as
every unoccupied square is on a line of the first type and every
occupied square is on a line of either the first or second type.  By
construction, the total degree of $f$ is $8k$.  If we group the
factors in $f$ according to slope, we see that $f$ can be rewritten as
\begin{equation}\label{eq:factor}
  f(x,y) = \prod_{j=1}^{2k} (x-\alpha_j)(y-\beta_j)(x-y-\gamma_j)(x+y-\delta_j)
\end{equation}
for suitable constants $\alpha_j,\beta_j,\gamma_j,\delta_j$.  From
equation~\eqref{eq:factor} and the binomial theorem, we conclude that
the coefficient of the top-degree term $x^{4k}y^{4k}$ is $\pm {2k
  \choose k}$, i.e., nonzero.

We now apply Theorem~\ref{theorem:CN} to $f(x,y)$, where $t_1=t_2=4k$
and $S_1=S_2=\{1,\ldots, 4k+1\}$, to obtain that there are $s_1 \in S_1, s_2 \in
S_2$ such that $f(s_1,s_2) \neq 0$.  We have reached a contradiction.
Therefore, the result holds when $n$ is congruent to $1$ modulo $4$.

Let us now consider $n=4k$, where $k$ is a positive integer.  We
proceed in a similar manner to the above.  Again, we consider a good
placement on $B_n$, this time of size $q \leq 4k-1$.  Similarly, let
$q'=4r+s$ denote the size of ${\cal Q}'$, where $r$ and $s$ are
integers with $0 \leq s \leq 3$.  As before, we define lines of the
first type and the second type, distributing those of the second type
as evenly as possible.  For each possible slope, the number of lines
is at most
\begin{equation}\label{equation:linenumber}
  g(r,s) = \left\lfloor \frac{4k-1-4r-s}{2} \right\rfloor + 
  \left\lceil \frac{4r+s}{4} \right\rceil.
\end{equation}
For $s \neq 1$, we have $g(r,s) \leq 2k-1$.  Likewise, for $s=1$ and
$r >0$, we have $g(r,s) \leq 2k-1$.  For these values, we proceed as
before, adding more lines so that there are $2k-1$ of each possible
slope.  We may now construct a polynomial of degree $8k-4$ and see
that the coefficient of the $x^{4k-1}y^{4k-3}$ term is nonzero.  As
before, applying Theorem~\ref{theorem:CN} we reach a contradiction
with $t_1=4k-1$ and $t_2=4k-3$ and $S_1=S_2=\{1, \ldots, 4k\}$.

We are left to consider the case where $s=1$ and $r=0$, i.e., $q'=1$.
In the process of defining lines of the
first type, since $q \leq 4k-1$ and $q'=1$ we may have $2k-1$ lines of each possible slope.  We
define one new line of the second type for the single queen in ${\cal
  Q}'$, giving it slope $\infty$.  Finally, we add more lines as
necessary so that there are precisely $2k$ with slope $\infty$ and
$2k-1$ for each of the other three slopes.  Our polynomial has degree
$8k-3$ and we consider the following leading term with nonzero
coefficient, ${2k-1 \choose k}x^{2k}y^{2k-1}(x^2)^{k-1}(-y^2)^k =
(-1)^k{2k-1 \choose k}x^{4k-2}y^{4k-1}$.  As before, applying
Theorem~\ref{theorem:CN} we reach a contradiction with $t_1=4k-2,
t_2=4k-1$ and $S_1=S_2=\{1, \ldots, 4k\}$.  This completes the proof in this case.

The cases of $n=4k+2$ and $n=4k+3$ follow in a similar manner and are left to the reader.  This
completes the proof.\qed

\section{A second proof.}\label{section:Harrisproof}

We now present a second proof to Theorem~\ref{theorem:main} for the
case $n$ is even and obtain a slightly weaker result when $n$ is
odd by showing that one needs at least $n-1$ queens; this proof is more elementary than the one given in Section
\ref{section:CNproof}.  While we arrived at it independently, many of
the ideas are to be found in a June 7, 1975 letter of John Harris to
Martin Gardner \cite{harris-letter}.  In the case when $n \equiv 3
\mod 4$, Harris only claimed $n-1$ queens are required.  A similar
proof for the no-two-in-a-line problem can be found in~\cite[Chapter
8]{watkins:04}.

We may refer to a square of $B_n$ by the coordinates $(x,y)$ of its
corresponding vertex.  

{\sc Proof:}  The claim is easily checked for $n=1$, so we assume $n\geq 2$.  Let ${\cal Q}$ be a good placement of size $q$ on $B_n$.  We distinguish between the lines of slope $0$ or
$\infty$ defined by ${\cal Q}$ and those of slope $\pm 1$.  To this
end, set $U\subseteq B_n$ 
to be the set of squares left uncovered by a
line of slope $0$ or $\infty$ and set ${\cal Q}''\subseteq {\cal Q}$ to
be those queens not involved in defining a line of slope $0$ or
$\infty$.  (Note that squares in $U$ may still be occupied by a queen
in ${\cal Q}''$.)  Write $q'' = |{\cal Q}''|$.  For any index $i \in \{1, \ldots, n\}$
(respectively $j \in \{1, \ldots, n\}$) let $C_i = \{(i,k)\in U :\, 1\leq k\leq n\}$
(respectively $R_j = \{(k,j)\in U :\, 1\leq k\leq n\}$).

\begin{figure}[htbp]
  \centering 
     {\scalebox{.15}{\includegraphics{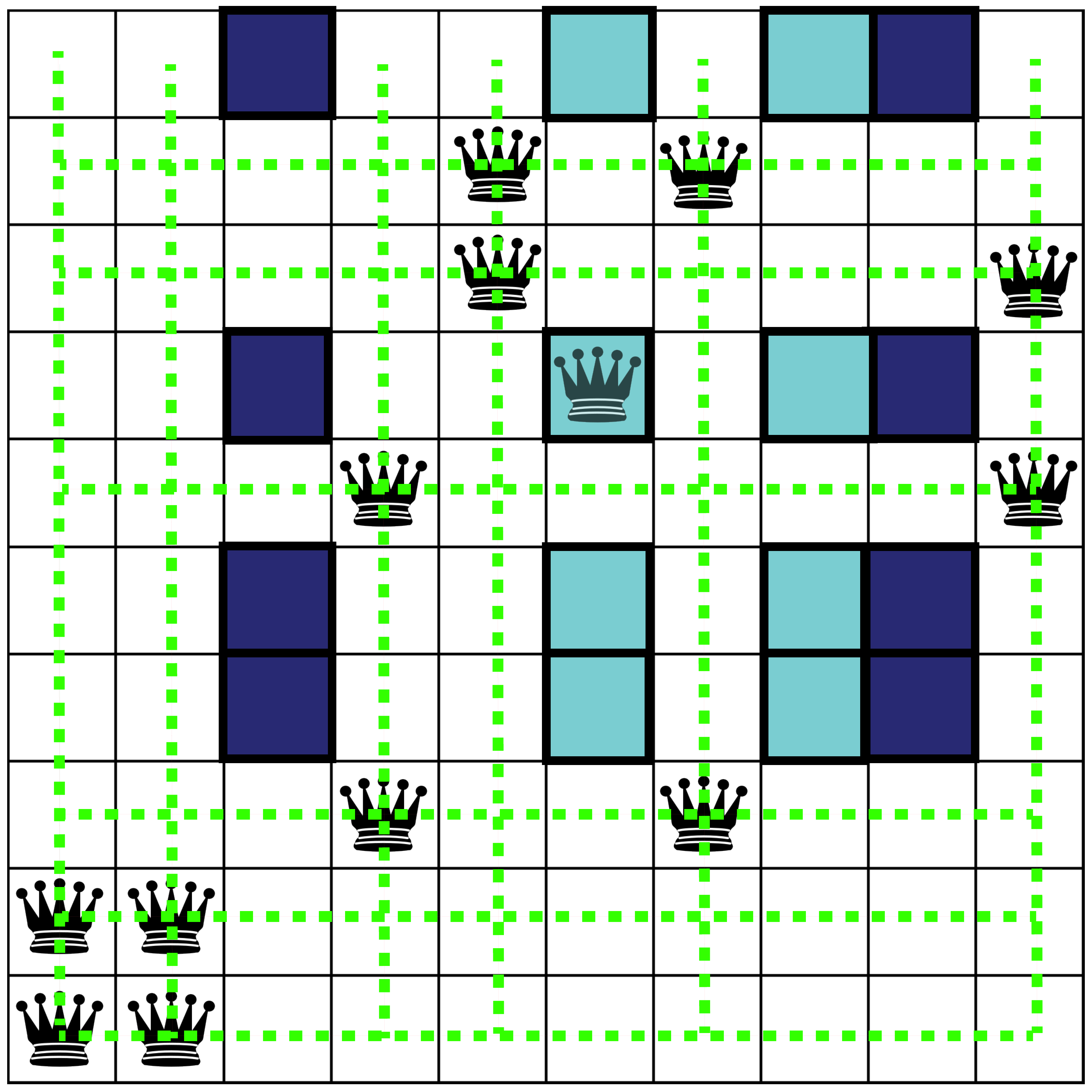}}}
     \caption{$q=13$, $q''=1$, $c=r=4$.  Squares of $U$ are shaded.  Dark shading indicates those that are also in $C_3 \cup
       C_9$.  The pale-shaded queen indicates the single queen in ${\cal
         Q}''$.}
      \label{fig:flyswat}
\end{figure}

The sets $C_i$ and $R_j$ keep track of the squares in $U$ for each
column and row.  Let $a < b$ be the minimum and maximum indices,
respectively, for which $C_i \neq \emptyset$.  Set $c$ to be the
number of the $C_i$ that are nonempty.  Define $a' < b'$ and $r$
analogously for the sets $R_j$.  Note that $c,r \geq n-
\frac{q-q''}{2}$.  In particular, $c\leq 1$ or $r\leq 1$ requires
$q\geq 2(n-1)$.  We therefore assume for the rest of the proof that
$r, c \geq 2$.  Without loss of generality, we may assume $b - a \geq
b' - a'$ as otherwise we may rotate the placement by $90^\circ$.
Figure~\ref{fig:flyswat} illustrates the various definitions of a $13$-queen good placement on a $10 \times 10$ chessboard.

As ${\cal Q}$ is good, the squares of $C_a \cup C_b$ are either
occupied or `attacked' via a pair of queens that would define a line
of slope $\pm 1$.  By definition, $|{\cal Q} \cap C_a| \leq 1, ~
|{\cal Q} \cap C_b| \leq 1$, and so $|{\cal Q} \cap (C_a \cup C_b)|
\leq \min\{q'', 2\}$.  There is at most one line of slope $+1$ that
attacks two squares of $C_a \cup C_b$ (the line would be a diagonal of
the `rectangle' formed by $C_a \cup C_b \cup R_{a'} \cup R_{b'}$).
Likewise, there is at most one line of slope $-1$ that attacks two
squares of $C_a \cup C_b$.  Each of the other lines of slope $\pm 1$
defined by ${\cal Q}$ attack at most one square of $C_a \cup C_b$.
The placement ${\cal Q}$ must therefore define at least
$2r-2-\min\{q'',2\}$ lines of slope $\pm 1$.  Furthermore,

\begin{equation}
  2r-2-\min\{q'',2\} \geq 2\left(n - \frac{q-q''}{2}\right)-2-q'' = 2n-q-2.
\end{equation}

Note that the $q$ queens of ${\cal Q}$ can define at most $q$ lines of
slope $\pm 1$.  Thus, $q \geq 2n-q-2$, and so $q \geq n-1$.

We now restrict $n$ to be even and we will reach a contradiction by
assuming that $q \leq n-1$.  As $n-1$ is odd, there are at most
$\frac{n-2}{2}$ lines of each possible slope defined by the placement
${\cal Q}$.  In particular, there are a total of at most $n-2$ lines
of slopes $\pm 1$.

If $q''=0$, then $r \geq n - \frac{n-2}{2} = \frac{n}{2}+1$, and so
$2r-2 \geq n$.  So, we need at least $n$ lines of slope $\pm 1$ --- a
contradiction.

If $q''>0$, then $r \geq n - \frac{q-q''}{2} \geq n- \frac{(n-1)-q''}{2}=
\frac{n}{2}+\frac{q''+1}{2}$.  We have $2r-2-\min\{q'',2\} \geq
2(\frac{n}{2}+\frac{q''+1}{2})-2-q'' = n-1$, and so we need at least
$n-1$ lines of slopes $\pm 1$ --- again, a contradiction.\qed

\section{Acknowledgments.}

The authors wish to thank Dan Archdeacon for some initial
conversations that led to this work.  We also wish to thank the Special Collections and
University Archives of Stanford University for helping us to access
Gardner's notes and correspondence.

The first author was supported by Middlebury College's Undergraduate Collaborative Research Fund.  The second author was partially supported by the National Science Foundation, Grant DMS-1100215.  The third author was partially supported by the National Science Foundation,  Grant DMS-0758057 and the National Security Agency, Grant H98230-10-1-0173.  The fourth author was partially supported by the National Security Agency, Grant H98230-09-1-0023 and the Simons Foundation (197419 to GSW).

\bigskip

\noindent\textbf{Alec S. Cooper} received his B.A. in mathematics from Middlebury College in May 2013.  He is particularly interested in algebra and other areas of discrete mathematics.

\noindent\textit{alecscooper@gmail.com}

\noindent\textbf{Oleg Pikhurko} received his Ph.D. in Mathematics from
Cambridge University in 2000.  He has an Erd\H{o}s number of two and
an Erd\H{o}s Lap number of two. Although he is the founder and CEO of
the Hedgehog Fund, it is unclear if he will have a finite Hedgehog Lap
number.

\noindent\textit{Mathematics Institute and DIMAP, University of Warwick,
Coventry, CV4 7AL, UK.\\
O.Pikhurko@warwick.ac.uk}

\noindent\textbf{John R. Schmitt} received his B.A. from Providence College in 1994, his M.S. from the University of Vermont in 1998, and his Ph.D. from Emory University in 2005.  He currently teaches at Middlebury College, where he devotes his research time to extremal combinatorics and graph theory.  He enjoys time spent with his wife and four children.

\noindent\textit{Department of Mathematics, Middlebury College, Middlebury, VT 05753. \\ jschmitt@middlebury.edu}

\noindent\textbf{Gregory S. Warrington}, an algebraic combinatorialist at the University of
Vermont, received his Ph.D. in Mathematics from Harvard University in
2001.  He likes to spend time with his family.

\noindent\textit{Department of Mathematics and Statistics, University of Vermont, Burlington, VT 05401. \\ gregory.warrington@uvm.edu}

\end{document}